\documentclass[11pt,letter]{article}%
\pdfoutput=1
\usepackage{amsmath}
\usepackage{amsfonts}
\usepackage{amssymb}
\usepackage{amsthm}
\usepackage{graphicx}%
\providecommand{\U}[1]{\protect\rule{.1in}{.1in}}
\newtheorem{theorem}{Theorem}[section]

\theoremstyle{definition}

\theoremstyle{remark}
\newtheorem{example}[theorem]{Example}
\newtheorem{remark}[theorem]{Remark}
\numberwithin{equation}{section}
\begin{document}

\title{Modulated electromagnetic fields in inhomogeneous media, hyperbolic pseudoanalytic functions and transmutations}
\author{Kira V. Khmelnytskaya$^{1}$, Vladislav V. Kravchenko$^{2}$, Sergii M.
Torba$^{2}$\thanks{The authors acknowledge the support from CONACYT, Mexico via the projects 166141 and 222478.}\\$^{1}${\footnotesize Faculty of Engineering, Autonomous University of
Queretaro, }\\{\footnotesize Cerro de las Campanas s/n, col. Las Campanas  Quer\'{e}taro,
Qro. C.P. 76010 M\'{e}xico}\\$^{2}${\footnotesize  Department of Mathematics, CINVESTAV del IPN, Unidad
Quer\'{e}taro }\\{\footnotesize Libramiento Norponiente \# 2000 Fracc. Real de Juriquilla
\ Quer\'{e}taro, Qro., CP 76230, M\'{e}xico}\\{\footnotesize vkravchenko@math.cinvestav.edu.mx}}
\maketitle

\begin{abstract}
The time-dependent Maxwell system describing electromagnetic wave propagation
in inhomogeneous isotropic media in the one-dimensional case reduces to a
Vekua-type equation for bicomplex-valued functions of a hyperbolic variable,
see \cite{KrRamirez2011}. Using this relation we solve the problem of the
transmission through an inhomogeneous layer of a normally incident
electromagnetic time-dependent plane wave. The solution is written in terms of
a pair of Darboux-associated transmutation operators \cite{KrT2012}, and
combined with the recent results on their construction \cite{KrT-CAOT},
\cite{KrT-Analytic} can be used for efficient computation of the transmitted
modulated signals. We develop the corresponding numerical method and
illustrate its performance with examples.

\end{abstract}

\section{Introduction}
In the present work a method for solving the problem of time-dependent
electromagnetic wave propagation through an isotropic inhomogeneous medium is
developed. Several ideas concerning such mathematical notions as bicomplex and
biquaternionic reformulations of electromagnetic models, hyperbolic Vekua-type
equations and transmutation operators from the theory of ordinary linear
differential equations are combined in our approach which results in a simple
and practical representation for solutions.

We observe that the 1+1 Maxwell system for inhomogeneous media can be
transformed into a hyperbolic Vekua equation. This gives us the possibility to
obtain the exact solution of the problem of a normally incident
electromagnetic time-dependent plane wave propagated through an inhomogeneous
layer in terms of a couple of Darboux-associated transmutation operators. This is a new representation for a solution of the classic problem.
Application of the recent results on the analytic approximation of such
operators allows us to write down the electromagnetic wave in an approximate
analytic form which is then used for numerical computation. The numerical
implementation of the proposed approach reduces to a certain recursive
integration and solution of an approximation problem, and can be based on the
usage of corresponding standard routines of such packages as Matlab.

In Section \ref{Sect2} we recall a biquaternionic reformulation of the Maxwell system
and use it to relate the 1+1 Maxwell system for inhomogeneous media with a
hyperbolic Vekua equation. We establish the equivalence between the
electromagnetic transmission problem and an initial-value problem for the
Vekua equation. In Section \ref{Sect3} we obtain the exact solution of the problem in
terms of a couple of the transmutation operators. In Section \ref{Sect4} the analytic
approximation of the exact solution is obtained in the case when the incident
wave is a partial sum of a trigonometric series. Other types of initial data
(and hence modulations) are discussed in Section \ref{Sect5}. Section \ref{Sect6} contains several
exactly solved examples used as test problems for the resulting numerical
method. In Section \ref{Sect7} we formulate some additional conclusions.

\section{The 1+1 Maxwell system as a Vekua equation and the problem statement}\label{Sect2}
The algebraic formalism in studying electromagnetic phenomena plays an
important role since Maxwell's original treatise in which Hamilton's
quaternions were present. In his PhD thesis of 1919 Lanczos wrote the Maxwell
system for a vacuum in the form of a single biquaternionic equation. This
elegant form of Maxwell' system was rediscovered in several posterior
publications, see, e.g., \cite{Imaeda}.

In \cite{Kr2003Progress} (see also \cite{KrAQA2003}) the Maxwell system
describing electromagnetic phenomena in isotropic inhomogeneous media was
written as a single biquaternionic equation. The Maxwell system has the form%
\begin{equation}
\operatorname{rot}\mathbf{H}=\varepsilon\partial_{t}\mathbf{E}+\mathbf{j,}
\label{Min1}%
\end{equation}%
\begin{equation}
\operatorname{rot}\mathbf{E}=-\mu\partial_{t}\mathbf{H}, \label{Min2}%
\end{equation}%
\begin{equation}
\operatorname{div}(\varepsilon\mathbf{E})=\rho, \label{Min3}%
\end{equation}%
\begin{equation}
\operatorname{div}\mathbf{(}\mu\mathbf{H})=0. \label{Min4}%
\end{equation}
\ Here $\varepsilon$ and $\mu$ are real-valued functions of spatial
coordinates, $\mathbf{E}$ and $\mathbf{H}$ are real-valued vector fields
depending on $t$ and spatial variables, the real-valued scalar function $\rho$
and vector function $\mathbf{j}$ characterize the distribution of sources of
the electromagnetic field. The following biquaternionic equation obtained in
\cite{Kr2003Progress} is equivalent to this system
\begin{equation}
\left( \frac{i}{c}\partial_{t}-D\right) \mathbf{V}+\mathbf{Vc}+\mathbf{V}%
^{\ast}\mathbf{Z}=\frac{\rho}{\sqrt{\varepsilon}}-i\sqrt{\mu}\mathbf{j.}
\label{Maxwell Quat}%
\end{equation}
Here
\[
\mathbf{c}=\frac{\operatorname*{grad}\sqrt{c}}{\sqrt{c}},\text{ \ }%
\mathbf{Z}=\frac{\operatorname*{grad}\sqrt{Z}}{\sqrt{Z}}\quad\text{and
}\mathbf{V}=\sqrt{\varepsilon}\mathbf{E}+i\sqrt{\mu}\mathbf{H}%
\]
where $c=\frac{1}{\sqrt{\varepsilon\mu}}$ is the wave propagation velocity and
$Z=\sqrt{\frac{\mu}{\varepsilon}}$ is the the intrinsic impedance of the
medium. All magnitudes in bold are understood as purely vectorial
biquaternions, and the asterisk denotes the complex conjugation (with respect
to the complex imaginary unit $i$). The operator $D$ is the main quaternionic
differential operator introduced by Hamilton himself and sometimes called the
Moisil-Theodoresco operator. It is defined on continuously differentiable
biquaternion-valued functions of the real variables $x_{1}$, $x_{2}$ and
$x_{3}$ according to the rule%
\[
Dq=\sum_{k=1}^{3}e_{k}\partial_{k}q,
\]
where $\partial_{k}=\frac{\partial}{\partial x_{k}}$ and $e_{k}$ are basic
quaternionic units.

In \cite{KrRamirez2011} with the aid of the representation of the Maxwell
system in the form (\ref{Maxwell Quat}) it was observed that in the sourceless
situation (i.e., $\rho$ and $\mathbf{j}$ are identically zeros) and when all
the magnitudes involved are independent of two spatial coordinates, say,
$x_{2}$ and $x_{3}$, $\varepsilon=\varepsilon(x_{1})$ and $\mu
=\operatorname*{Const}$, the Maxwell system is equivalent to the following
Vekua-type equation
\begin{equation}
\partial_{\overline{z}}W-\frac{f^{\prime}}{2f}\overline{W}=0 \label{BicompW}%
\end{equation}
where $\partial_{\overline{z}}=\frac{1}{2}(\partial_{\xi}-j\partial_{t})$, $j$
is a hyperbolic imaginary unit, $j^{2}=1$ commuting with $i$, $W$ is a
bicomplex-valued function of the real variables $\xi$ and $t$, $W=u+vj$ and
$u$, $v$ are complex valued (containing the imaginary unit $i$). The function
$f$ is real valued and depends on $\xi$ only. The conjugation with respect to
$j$ is denoted by the bar, $\overline{W}=u-vj$.

The Maxwell system in this case can be written in the form%
\begin{equation}
\varepsilon(x)\partial_{t}\mathcal{E}=i\partial_{x}\mathcal{H},\quad
i\partial_{x}\mathcal{E}=-\mu\partial_{t}\mathcal{H} \label{Maxwell1+1}%
\end{equation}
where $\mathcal{E}=E_{2}+iE_{3}$, $\mathcal{H}=H_{2}+iH_{3}$, $x=x_{1}$. The
relation between (\ref{Maxwell1+1}) and (\ref{BicompW}) involves the change of
the independent variable $\xi(x)=\sqrt{\mu}\int_{0}^{x}\sqrt{\varepsilon
(s)}ds$. The function $f$ in (\ref{BicompW}) is related to $\varepsilon$ and
$\mu$ by the equality $f(\xi)=\sqrt{\widetilde{c}(0)}/\sqrt{\widetilde{c}%
(\xi)}$ where and below the tilde means that a function of $x$ is written as a
function of $\xi$, $\widetilde{c}(\xi(x))=c(x)$. The function $W$ is written
in terms of $E$ and $H$ as follows%
\begin{equation}
W(\xi,t)=\sqrt{\widetilde{c}(\xi)}\left(  \sqrt{\widetilde{\varepsilon}(\xi
)}\widetilde{\mathcal{E}}(\xi,t)+ij\sqrt{\mu}\widetilde{\mathcal{H}}\right)  .
\label{Relation W}%
\end{equation}

Let us consider the problem of a normally incident plane wave transmission
through an inhomogeneous medium (see, e.g., \cite[Chapter 8]{OstrovskyPotapov}%
). The electromagnetic field $\mathcal{E}$ and $\mathcal{H}$ is supposed to be
known at $x=0$,%
\begin{equation}
\mathcal{E}(0,t)=\mathcal{E}_{0}(t)\quad\text{and}\quad\mathcal{H}%
(0,t)=\mathcal{H}_{0}(t)\text{,\quad}t\in\lbrack\alpha,\beta].
\label{init cond}%
\end{equation}
We assume $\mathcal{E}_{0}$ and $\mathcal{H}_{0}$ to be continuously
differentiable functions.

The problem (\ref{Maxwell1+1}), (\ref{init cond}) can be reformulated in terms
of the function (\ref{Relation W}). Find a solution of (\ref{BicompW})
satisfying the condition%
\begin{equation}
W(0,t)=W_{0}(t) \label{initial condition}%
\end{equation}
where
\begin{equation}
W_{0}=\sqrt{c(0)\varepsilon(0)}\mathcal{E}_{0}+ij\sqrt{c(0)\mu}\mathcal{H}_{0}
\label{W0}%
\end{equation}
is a given continuously differentiable function.

\section{The representation of the solution}\label{Sect3}

First, let us consider the elementary problem%
\begin{align}
w_{\overline{z}} & =0,\label{CR hyper}\\
w(0,t) & =w_{0}(t). \label{initial condition w}%
\end{align}
The hyperbolic Cauchy-Riemann system (\ref{CR hyper}) was studied in several
publications (see, e.g., \cite{Lavrentyev and Shabat}, \cite{MotterRosa},
\cite{Wen} and more recent \cite{KrT-CAOT}). Its general solution can be
written in the form%
\[
w(\xi,t)=P^{+}\Phi(t+\xi)+P^{-}\Psi(t-\xi)
\]
where $\Phi$ and $\Psi$ are arbitrary continuously differentiable scalar
functions, $P^{\pm}=\frac{1}{2}\left(  1\pm j\right)  $.

For the scalar components of $w$ we introduce the notations
\[
\mathcal{R}(w)=u=\frac{1}{2}(w+\overline{w})\quad\text{and\quad}%
\mathcal{I}(w)=v=\frac{1}{2j}(w-\overline{w}).
\]

When $\xi=0$ we obtain $w(0,t)=P^{+}\Phi(t)+P^{-}\Psi(t)$. Hence the unique
solution of the Cauchy problem (\ref{CR hyper}), (\ref{initial condition w})
has the form%
\begin{equation}%
\begin{split}
w(\xi,t)  &  =P^{+}w_{0}^{+}(t+\xi)+P^{-}w_{0}^{-}(t-\xi)\\
&  =\frac{1}{2}\left(  w_{0}^{+}(t+\xi)+w_{0}^{-}(t-\xi)+j\left(  w_{0}%
^{+}(t+\xi)-w_{0}^{-}(t-\xi)\right)  \right)
\end{split}
\label{elementary sol}%
\end{equation}
where
\begin{equation}
w_{0}^{\pm}:=\mathcal{R}\left(  w_{0}\right)  \pm\mathcal{I}(w_{0}).
\label{W+-}%
\end{equation}

In \cite{KrT-CAOT} there was established a relation between solutions of
(\ref{CR hyper}) and solutions of (\ref{BicompW}). Any solution of
(\ref{BicompW}) can be represented in the form
\begin{equation}
W(\xi,t)=T_{f}\left[  \mathcal{R}\left(  w(\xi,t)\right)  \right]
+jT_{1/f}\left[  \mathcal{I}\left(  w(\xi,t)\right)  \right]
\label{W transmut w}%
\end{equation}
where $w$ is a solution of (\ref{CR hyper}), $T_{f}$ and $T_{1/f}$ are
Darboux-associated transmutation operators defined in \cite{KrT2012}, see also
\cite{CKM2012} and \cite{KrT-CAOT}, and applied with respect to the variable $\xi$.
Both operators have the form of Volterra
integral operators,
\[
T_{f}u(\xi)=u(\xi)+\int_{-\xi}^{\xi}\mathbf{K}_{f}(\xi,\tau)u(\tau)d\tau
\]
with continuously differentiable kernels. Moreover, the operators $T_{f}$ and
$T_{1/f}$ preserve the value at $\xi=0$ giving additionally to
(\ref{W transmut w}) the relation $W(0,t)=w(0,t)$. This together with
(\ref{elementary sol}) allows us to write down the unique solution of the
problem (\ref{BicompW}), (\ref{initial condition}) in the form
\begin{equation}
W(\xi,t)=\frac{1}{2}T_{f}\left[  W_{0}^{+}(t+\xi)+W_{0}^{-}(t-\xi)\right]
+\frac{j}{2}T_{1/f}\left[  W_{0}^{+}(t+\xi)-W_{0}^{-}(t-\xi)\right]  .
\label{Solution initial problem}%
\end{equation}
In what follows we use the convenience of this representation and the recent
results \cite{KrT-CAOT}, \cite{KrT-Analytic} on the construction of the
operators $T_{f}$ and $T_{1/f}$.

\begin{remark}
The single-wave approximation of the solution of (\ref{BicompW}),
(\ref{initial condition}) (see \cite[Subsection 8.5.2]{OstrovskyPotapov}) can
be obtained from the representation (\ref{Solution initial problem}) by
removing integrals from the definition of $T_{f}$ and $T_{1/f}$ or in other
words, replacing $T_{f}$ and $T_{1/f}$ by the identity operator.
\end{remark}

\section{Approximation of modulated waves, the simplest initial data}\label{Sect4}
Initial data interesting in practical problems correspond to
modulated electromagnetic waves which are represented as partial sums of
trigonometric series (other types of initial data are discussed in the next section). In other words, consider initial data of
the form
\begin{equation}
\mathcal{E}_{0}(t)=\sum_{m=-M}^{M}\alpha_{m}e^{i(\omega_{0}+m\omega)t}%
\quad\text{and}\quad\mathcal{H}_{0}(t)=\sum_{m=-M}^{M}\beta_{m}e^{i(\omega
_{0}+m\omega)t}. \label{EHmodulated}%
\end{equation}
This leads to a similar form for the initial data $W_{0}$ in
(\ref{initial condition}),
\begin{equation}
W_{0}(t)=\sum_{m=-M}^{M}\gamma_{m}e^{i(\omega_{0}+m\omega)t} \label{W0exp}%
\end{equation}
where the bicomplex numbers $\gamma_{m}$ are related to $\alpha_{m}$,
$\beta_{m}\in\mathbb{C}$ as follows%
\[
\gamma_{m}=\sqrt{c(0)}\left(  \sqrt{\varepsilon(0)}\alpha_{m}+ij\sqrt{\mu
}\beta_{m}\right)  .
\]
Due to (\ref{Solution initial problem}), the unique solution of the problem
(\ref{BicompW}), (\ref{initial condition}) with $W_{0}$ given by (\ref{W0exp})
can be written in the form%
\begin{equation}
\begin{split}
W(\xi,t)  &  =\frac{1}{2}\left(\sum_{m=-M}^{M}e^{i(\omega_{0}+m\omega)t}\left(
\gamma_{m}^{+}T_{f}\left[  e^{i(\omega_{0}+m\omega)\xi}\right]  +\gamma
_{m}^{-}T_{f}\left[  e^{-i(\omega_{0}+m\omega)\xi}\right]  \right)\right. \\
&  +j\left.\sum_{m=-M}^{M}e^{i(\omega_{0}+m\omega)t}\left(  \gamma_{m}^{+}%
T_{1/f}\left[  e^{i(\omega_{0}+m\omega)\xi}\right]  -\gamma_{m}^{-}%
T_{1/f}\left[  e^{-i(\omega_{0}+m\omega)\xi}\right]  \right)  \right)
\end{split}\label{Wxit}
\end{equation}
where $\gamma_{m}^{\pm}=\sqrt{c(0)}\left(  \sqrt{\varepsilon(0)}\alpha_{m}\pm
i\sqrt{\mu}\beta_{m}\right)  $.

Although the explicit form of the operators $T_{f}$ and $T_{1/f}$ is usually
unknown, in \cite{KrT-Analytic} it was shown how their kernels can be
approximated by means of generalized wave polynomials. In particular, for the
images of the functions $e^{\pm i(\omega_{0}+m\omega)\xi}$ the following
approximate representations are valid%
\begin{equation}%
\begin{split}
T_{f}\left[  e^{\pm i(\omega_{0}+m\omega)\xi}\right]   &  \cong e^{\pm
i(\omega_{0}+m\omega)\xi}+2\sum_{n=0}^{N}a_{n}\sum_{\text{even }k=0}^{n}%
\binom{n}{k}\varphi_{n-k}(\xi)\int_{0}^{\xi}\tau^{k}\cos(\omega_{0}%
+m\omega)\tau\,d\tau\\
& \quad\pm2i\sum_{n=1}^{N}b_{n}\sum_{\text{odd }k=1}^{n}\binom{n}{k}%
\varphi_{n-k}(\xi)\int_{0}^{\xi}\tau^{k}\sin(\omega_{0}+m\omega)\tau\,d\tau
\end{split}
\label{TfExp}%
\end{equation}
and%
\begin{equation}%
\begin{split}
T_{1/f}\left[  e^{\pm i(\omega_{0}+m\omega)\xi}\right]   &  \cong e^{\pm
i(\omega_{0}+m\omega)\xi}-2\sum_{n=0}^{N}b_{n}\sum_{\text{even }k=0}^{n}%
\binom{n}{k}\psi_{n-k}(\xi)\int_{0}^{\xi}\tau^{k}\cos(\omega_{0}+m\omega
)\tau\,d\tau\\
& \quad\mp2i\sum_{n=1}^{N}a_{n}\sum_{\text{odd }k=1}^{n}\binom{n}{k}\psi
_{n-k}(\xi)\int_{0}^{\xi}\tau^{k}\sin(\omega_{0}+m\omega)\tau\,d\tau.
\end{split}
\label{T1/fExp}%
\end{equation}
Here all integrals are easily calculated in a closed form, the functions
$\varphi_{n}$ and $\psi_{n}$ are defined as follows. Consider two sequences of
recursive integrals%
\begin{equation}
X^{(0)}(x)\equiv1,\qquad X^{(n)}(x)=n\int_{0}^{x}X^{(n-1)}(s)\left(
f^{2}(s)\right)  ^{(-1)^{n}}\,\mathrm{d}s,\qquad n=1,2,\ldots\label{Xn}%
\end{equation}
and
\begin{equation}
\widetilde{X}^{(0)}\equiv1,\qquad\widetilde{X}^{(n)}(x)=n\int_{0}%
^{x}\widetilde{X}^{(n-1)}(s)\left(  f^{2}(s)\right)  ^{(-1)^{n-1}}%
\,\mathrm{d}s,\qquad n=1,2,\ldots. \label{Xtiln}%
\end{equation}

The two families of functions $\left\{  \varphi_{k}\right\}  _{k=0}^{\infty}$
and $\left\{  \psi_{k}\right\}  _{k=0}^{\infty}$ are constructed according to
the rules
\begin{equation}
\varphi_{k}(x)=%
\begin{cases}
f(x)X^{(k)}(x), & k\text{\ odd},\\
f(x)\widetilde{X}^{(k)}(x), & k\text{\ even},
\end{cases}
\label{phik}%
\end{equation}
and
\begin{equation}
\psi_{k}(x)=%
\begin{cases}
\dfrac{\widetilde{X}^{(k)}(x)}{f(x)}, & k\text{\ odd,}\\
\dfrac{X^{(k)}(x)}{f(x)}, & k\text{\ even}.
\end{cases}
\label{psik}%
\end{equation}

Finally, the coefficients $a_{n}$ and $b_{n}$ are obtained by solving an
approximation problem described in \cite{KrT-Analytic}.

As an important feature of the representations (\ref{TfExp}) and
(\ref{T1/fExp}) in \cite{KrT-Analytic} there were obtained estimates of their
accuracy uniform with respect to the parameter $(\omega_{0}+m\omega)$.

Note that the direct evaluation of expressions \eqref{TfExp} and \eqref{T1/fExp} requires $O(N^2)$ algebraic operations for each pair of $m$ and $\xi$ leading to the computation complexity $O(MN^2)$ for the evaluation of \eqref{Wxit}. The change of summation order in \eqref{TfExp} and \eqref{T1/fExp} allows one to evaluate \eqref{Wxit} with the computation complexity of $O(NM+N^2)$. For example, consider
\begin{align*}
&\sum_{n=0}^N a_n  \sum_{\text{even }k=0}^{n}%
\binom{n}{k}\varphi_{n-k}(\xi)\int_{0}^{\xi}\tau^{k}\cos(\omega_{0}%
+m\omega)\tau\,d\tau\\
&\quad =\sum_{\text{even }k=0}^{N}\left(\sum_{n=k}^N a_n\binom nk \varphi_{n-k}(\xi)\right)\int_0^\xi \tau^k\cos(\omega_0+m\omega)\tau\,d\tau.
\end{align*}
The coefficients $\sum_{n=k}^N a_n\binom nk \varphi_{n-k}(\xi)$ can be precomputed once for every $\xi$ in $O(N^2)$ operations and the outer sum requires $O(N)$ operations.

\section{Other types of initial data}\label{Sect5}
The proposed method is not restricted exclusively to initial data which can be represented or closely approximated by \eqref{EHmodulated}. When one can efficiently calculate (at least numerically) the indefinite integrals
\[
\int \tau^k W_0^+(\tau)\,d\tau\qquad \text{and}\qquad \int \tau^k W_0^-(\tau)\,d\tau,\quad k=0,\ldots,N,
\]
the approximations of the transmutation operators can be used. One of the examples of such initial data arises in digital signal transmission. For widely used modulations like phase-shift keying or QAM, the transmitted signal can be represented as
\[
s(t) = \sum_{n=0}^{M}\left[c_n\cos \omega_0 t+ s_n\sin\omega_0 t\right]\mathbf{1}_{[nf_s ,(n+1)f_s)}(t),
\]
where $\mathbf{1}_{[a,b)}$ is the characteristic function of the interval $[a,b)$, $\omega_0$ is the carrier frequency, $f_s$ is the symbol rate and the coefficients $s_n$ and $c_n$ encode transmitted information. Other examples include Gaussian RF pulses and linear frequency modulation (also known as chirp modulation) used for radars. We refer the reader to \cite{Shmaliy} for further details.

Returning to \eqref{Solution initial problem}, consider $T_f[W_0^+(t+\xi)]$. In \cite{KrT-Analytic} the following approximation was constructed
\begin{equation}
    \begin{split}
       T_f[W_0^+(t+\xi)](\xi) &\cong W_0^+(t+\xi)+\sum_{n=0}^N a_n \sum_{\text{even } k=0}^n \binom nk \varphi_{n-k}(\xi)\int_{-\xi}^{\xi} \tau^k W_0^+(t+\tau)\, d\tau \\
         & \quad + \sum_{n=1}^N b_n \sum_{\text{odd } k=1}^n \binom nk \varphi_{n-k}(\xi)\int_{-\xi}^{\xi} \tau^k W_0^+(t+\tau)\, d\tau.
     \end{split}\label{TfW}
\end{equation}
We have
\[
\int_{-\xi}^{\xi} \tau^k W_0^+(t+\tau)\, d\tau = \int_{t-\xi}^{t+\xi}(z-t)^k W_0^+(z)\,dz = \sum_{\ell=0}^k \binom k\ell (-1)^{k-\ell} t^{k-\ell}\int_{t-\xi}^{t+\xi} z^\ell W_0^+(z)\,dz,
\]
hence reordering terms in \eqref{TfW} one obtains
\begin{equation}\label{TfWfin}
    \begin{split}
      T_f[W_0^+(t+\xi)] & \cong  W_0^+(t+\xi)+\sum_{\ell=0}^N \left(\int_{t-\xi}^{t+\xi} z^\ell W_0^+ (z)\, dz\right)\\
       &\quad \times \sum_{n=\ell}^N \Bigg\{a_n\sum_{\text{even }k\ge \ell}^n +b_n\sum_{\text{odd }k\ge \ell}^n\Bigg\} (-1)^{k-\ell}\binom nk \binom k\ell \varphi_{n-k}(\xi) t^{k-\ell},
    \end{split}
\end{equation}
an expression which can be easily evaluated once the indefinite integrals of $z^\ell W_0^+ (z)$ are known. Similarly,
\begin{equation}\label{TfWmfin}
    \begin{split}
      T_f[W_0^-(t-\xi)] & \cong  W_0^-(t-\xi)+\sum_{\ell=0}^N \left(\int_{t-\xi}^{t+\xi} z^\ell W_0^- (z)\, dz\right)\\
       &\quad \times \sum_{n=\ell}^N \Bigg\{a_n\sum_{\text{even }k\ge \ell}^n +b_n\sum_{\text{odd }k\ge \ell}^n\Bigg\} (-1)^{\ell}\binom nk \binom k\ell \varphi_{n-k}(\xi) t^{k-\ell}.
    \end{split}
\end{equation}
To obtain approximations for $T_{1/f}[W_0^\pm(t\pm\xi)]$ on can change in \eqref{TfWfin} and \eqref{TfWmfin} the coefficients $a_n$ and $b_n$ by $-b_n$ and $-a_n$ respectively and change functions $\phi_n$ by $\psi_n$.

\section{Examples and numerical tests}\label{Sect6}

\begin{example}
\label{Ex1} Let us consider the system (\ref{Maxwell1+1}) with the
permittivity of the form 
\begin{equation}\label{ExEpsilon}
\varepsilon(x)=(\alpha x+\beta)^{-2},
\end{equation}
where $\alpha$
and $\beta$ are some real numbers, such that $\alpha x+\beta$ does not vanish
on the interval of interest and $\varepsilon>0$. Then $\xi=\sqrt{\mu}\int
_{0}^{x}\sqrt{\varepsilon(s)}ds=\frac{\sqrt{\mu}}{\alpha}\log\frac{\alpha
x+\beta}{\beta}$. Hence
\[
x=\frac{\beta}{\alpha}\left(  e^{\frac{\alpha\xi}{\sqrt{\mu}}}-1\right)
\quad\text{and}\quad\widetilde{\varepsilon}(\xi)=\frac{1}{\beta^{2}}%
e^{-\frac{2\alpha\xi}{\sqrt{\mu}}},\quad\widetilde{c}(\xi)=\frac{\beta}%
{\sqrt{\mu}}e^{\frac{\alpha\xi}{\sqrt{\mu}}},\quad f(\xi)=e^{-\frac{\alpha\xi
}{2\sqrt{\mu}}}.
\]
In this case the Vekua equation (\ref{BicompW}) has the form
\begin{equation}
\partial_{\overline{z}}W+\gamma\overline{W}=0 \label{VekuaConst}%
\end{equation}
where the coefficient $\gamma$ is constant, $\gamma=\alpha/\left(  4\sqrt{\mu
}\right)  $.

The Vekua equation (\ref{BicompW}) is a special case of the main Vekua
equation $\partial_{\overline{z}}W-\frac{f_{\overline{z}}}{f}\overline{W}=0$
(see \cite{KRT}, \cite{APFT}). For its solutions one has that $\mathcal{R}(W)$
satisfies the equation $\left(  \square-q\right)  u=0$ where $\square
:=\partial_{\xi}^{2}-\partial_{t}^{2}$ and $q=\square f/f$, and $\mathcal{I}%
(W)$ satisfies the equation $\left(  \square-r\right)  v=0$ where
$r=\square(f^{-1})/(f^{-1})$. The relation between $\mathcal{R}(W)$ and
$\mathcal{I}(W)$ is akin to the relation between harmonic conjugate functions
and can be found in \cite{KRT}, \cite{APFT}.

In the case of equation (\ref{VekuaConst}) $q=r=\alpha^{2}/(4\mu)$. Consider a
solution $u$ of the equation $\left(  \square-\frac{\alpha^{2}}{4\mu}\right)
u=0$ in the form $u(\xi,t)=e^{2i\gamma t}(A\xi+B)$ where $A$ and $B$ are
arbitrary real constants. Then a solution $W$ of (\ref{VekuaConst}) such that
$\mathcal{R}(W)=u$ can be chosen in the form (see \cite{KRT}, \cite{APFT})
\begin{equation}
W(\xi,t)=e^{2i\gamma t}\left(  \left(  1-ij\right)  (A\xi+B)-ij\frac
{A}{2\gamma}\right)  .\label{solVekuaConst}%
\end{equation}
This leads to the following solution of the Maxwell system
\[
\widetilde{\mathcal{E}}(\xi,t)=\frac{1}{\sqrt{\widetilde{c}(\xi)\widetilde
{\varepsilon}(\xi)}}\mathcal{R}(W(\xi,t))=\frac{e^{2i\gamma t}}{\sqrt
{\widetilde{c}(\xi)\widetilde{\varepsilon}(\xi)}}(A\xi+B)
\]
and%
\[
\widetilde{\mathcal{H}}(\xi,t)=-\frac{i}{\sqrt{\widetilde{c}(\xi)\mu}%
}\mathcal{I}(W(\xi,t))=-\frac{e^{2i\gamma t}}{\sqrt{\widetilde{c}(\xi)\mu}%
}\left(  A\xi+B+\frac{A}{2\gamma}\right)
\]
or in terms of the variables $x$ and $t$,
\begin{equation}\label{Ex1Ext}
\mathcal{E}(x,t)=\sqrt[4]{\mu}\sqrt{\alpha x+\beta}e^{\frac{i\alpha t}%
{2\sqrt{\mu}}}\left(  \frac{\sqrt{\mu}}{\alpha}A\log\frac{\alpha x+\beta
}{\beta}+B\right)
\end{equation}
and%
\begin{equation}\label{Ex1Hxt}
\mathcal{H}(x,t)=-\frac{e^{\frac{i\alpha t}{2\sqrt{\mu}}}}{\sqrt[4]{\mu}%
\sqrt{\alpha x+\beta}}\left(  \frac{\sqrt{\mu}}{\alpha}A\log\frac{\alpha
x+\beta}{\beta}+B+\frac{2A\sqrt{\mu}}{\alpha}\right)  .
\end{equation}
Thus the functions \eqref{Ex1Ext} and \eqref{Ex1Hxt} are the exact solutions of the Maxwell system \eqref{Maxwell1+1} with $\varepsilon(x)$ given by \eqref{ExEpsilon} and with the initial conditions
\[
\mathcal{E}(0,t)=\sqrt[4]{\mu}\sqrt{\beta}e^{\frac{i\alpha t}%
{2\sqrt{\mu}}}B\quad \text{and}\quad \mathcal{H}(0,t)=
-\frac{e^{\frac{i\alpha t}{2\sqrt{\mu}}}}{\sqrt[4]{\mu}
\sqrt{\beta}}\left(  B+\frac{2A\sqrt{\mu}}{\alpha}\right).
\]
Below we present the results of numerical solution of the same system by the proposed method.

For the numerical experiment we considered an interval $[0,5]$ for both $x$
and $t$ and took $\alpha=2$, $\beta=1$, $\mu=1$ and $A=1$, $B=3$. The initial
condition in this example corresponds to the case $M=0$ and $\omega
_{0}=2\gamma$ in (\ref{W0exp}). The permittivity $\varepsilon(x)$ was
approximated on uniform mesh of $5001$ points. The new variable $\xi$ was
obtained by the modified 6 point Newton-Cottes integration formula, see
\cite{CKT2013} for details. The same integration formula was used for
calculation of all other integrals, i.e., for computation of \eqref{Xn} and
\eqref{Xtiln} and for the integration of the potential $q$ required to obtain
approximations \eqref{TfExp} and \eqref{T1/fExp}. Note that the integration
with respect to the variable $\xi$ requires integration over a non-uniform
mesh, however such inconvenience can be easily avoided observing that
$\int_{0}^{\xi(x)} \widetilde g(\xi)\,d\xi= \int_{0}^{x} g(s) \xi^{\prime
}(s)\,ds = \int_{0}^{x} g(s)\sqrt{\mu\varepsilon(s)}\,ds$ for any function
$g(x)=\widetilde g(\xi(x))$.

All calculations were performed using Matlab 2012 in the machine precision.
The exact expressions were used only for the function $\varepsilon(x)$ and its
derivatives, all other functions involved were computed numerically.

The developed program found the optimal value of $N$ for the approximations
\eqref{TfExp} and \eqref{T1/fExp} to be $N=14$. The computation time required
was 0.4 seconds. On Figure \ref{fig:Ex1} we show the graphs of the absolute
errors of the computed $\mathcal{E}(x,t)$ and $\mathcal{H}(x,t)$.
\begin{figure}[tbh]
\centering
\includegraphics[
natwidth=900,
natheight=720,
width=3in,
height=2.4in]{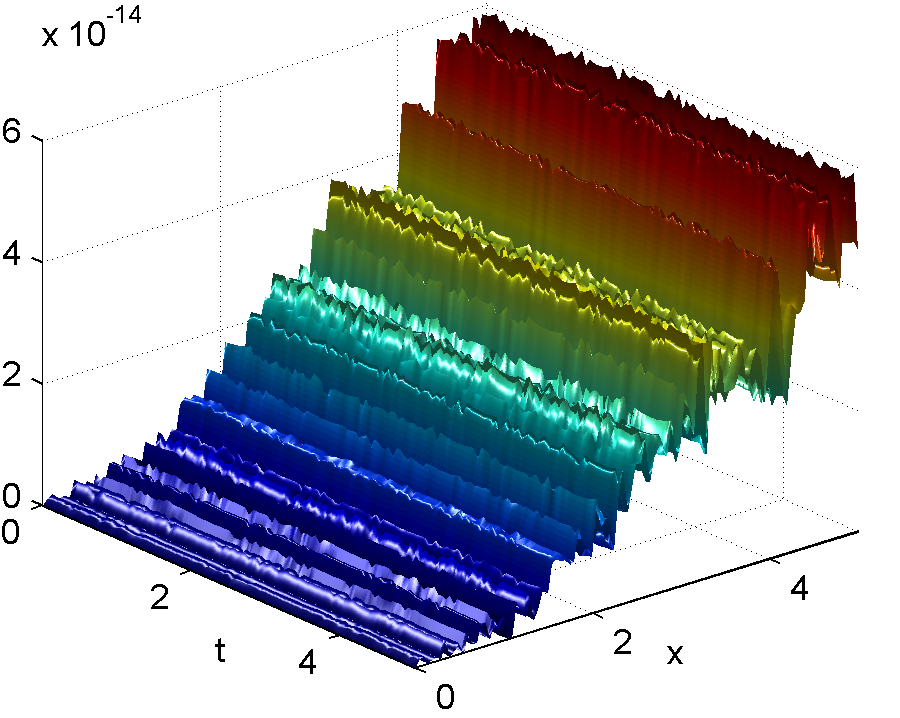} \quad\includegraphics[
natwidth=900,
natheight=718,
width=3in,
height=2.4in]{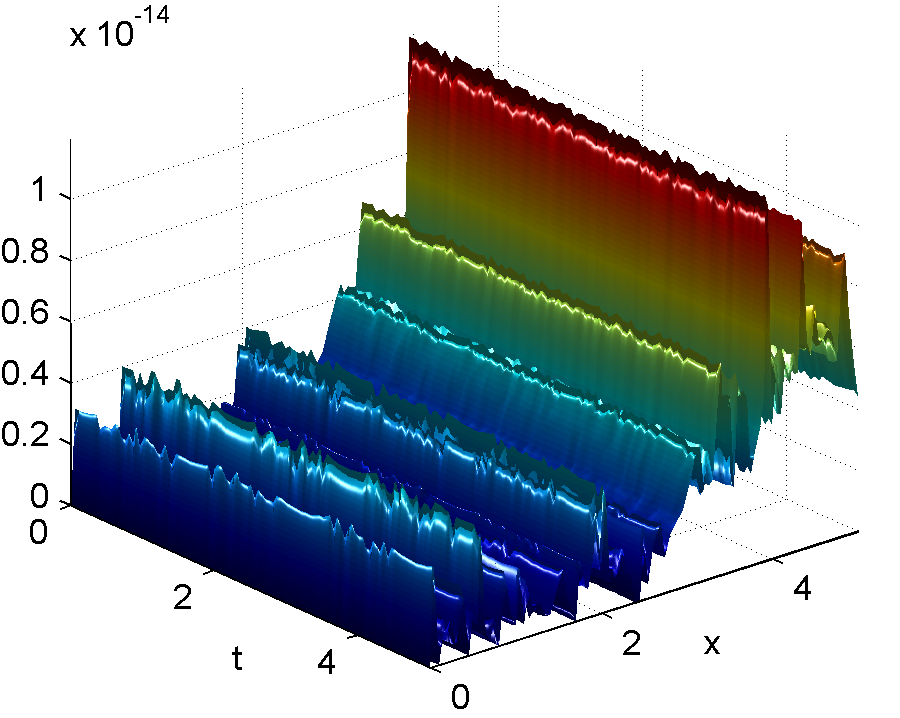}\caption{Graphs of the absolute errors of
$\mathcal{E}(x,t)$ (on the left) and $\mathcal{H}(x,t)$ (on the right) from
Example \ref{Ex1}.}%
\label{fig:Ex1}%
\end{figure}
\end{example}

\begin{example}\label{Ex2}
Let us consider the same parameters of the medium as in Example \ref{Ex1} but
instead of choosing a solution of (\ref{VekuaConst}) in the form
(\ref{solVekuaConst}) we take the solution%
\begin{equation}\label{Ex2Wxit}
W(\xi,t)=Ae^{i\Omega t}\left(  e^{D\xi}+\frac{D+C}{D-C}e^{-D\xi}%
+\frac{2ij\Omega}{D-C}\sinh D\xi\right)  .
\end{equation}
Here $A$ and $\Omega$ are arbitrary constants, $C=\frac{\alpha}{2\sqrt{\mu}}$,
$D=i\sqrt{\Omega^{2}-C^{2}}$. It is obtained similarly to (\ref{solVekuaConst}%
) starting with a solution $u$ of the equation $\left(  \square-\frac
{\alpha^{2}}{4\mu}\right)  u=0$ in the form $u(\xi,t)=Ae^{i\Omega t}(e^{D\xi
}+\frac{D+C}{D-C}e^{-D\xi})$.

We have then that
\[
\mathcal{E}(x,t)=A\sqrt[4]{\mu}\sqrt{\alpha x+\beta}e^{i\Omega t}\left(
\left(  \frac{\alpha x+\beta}{\beta}\right)  ^{\frac{D\sqrt{\mu}}{\alpha}%
}+\frac{D+C}{D-C}\left(  \frac{\alpha x+\beta}{\beta}\right)  ^{-\frac
{D\sqrt{\mu}}{\alpha}}\right)
\]
and%
\[
\mathcal{H}(x,t)=\frac{A}{D-C}\frac{\Omega\,e^{i\Omega t}}{\sqrt[4]{\mu}%
\sqrt{\alpha x+\beta}}\left(  \left(  \frac{\alpha x+\beta}{\beta}\right)
^{\frac{D\sqrt{\mu}}{\alpha}}-\left(  \frac{\alpha x+\beta}{\beta}\right)
^{-\frac{D\sqrt{\mu}}{\alpha}}\right)  
\]
satisfy the Maxwell system \eqref{Maxwell1+1} with the permittivity \eqref{ExEpsilon} and the initial conditions
\[
\mathcal{E}(0,t)=\frac{2AD}{D-C}\sqrt[4]{\mu}\sqrt{\beta}e^{i\Omega t}\quad\text{and}\quad 
\mathcal{H}(0,t)=0.
\]

For the numerical calculation we considered the same values of the parameters $\alpha$, $\beta$, $\mu$ as in Example \ref{Ex1} and took the interval $[0,6]$ for both $x$ and $t$. For the initial condition we took the sum of four terms, each of the form \eqref{Ex2Wxit} having $\Omega_1=-\Omega_2=C+1$, $\Omega_3=-\Omega_4=C+2$. Since the expression \eqref{Ex2Wxit} for $\xi=0$ reduces to $W(0,t)=\frac{2AD}{D-C}e^{i\Omega t}$, we took $A_i = \frac{D_i-C}{D_i}$, $i=1,\ldots,4$ and obtained initial conditions $W^\pm_0(t)=4 \cos (C+1)t + 4\cos (C+2)t$.

For this example the optimal $N$ was equal to 13 and the whole computation time was $0.3$ seconds. On Figure \ref{fig:Ex2} we present the graphs of the initial condition and the computed $\mathcal{E}(x,t)$. The absolute errors of the computed $\mathcal{E}(x,t)$ and $\mathcal{H}(x,t)$ were less than $1.1\cdot 10^{-13}$ and $9\cdot 10^{-15}$ respectively.
\begin{figure}[tbh]
\centering
\includegraphics[
width=2.2in,
height=2.2in,
bb = 216   306   396   486]{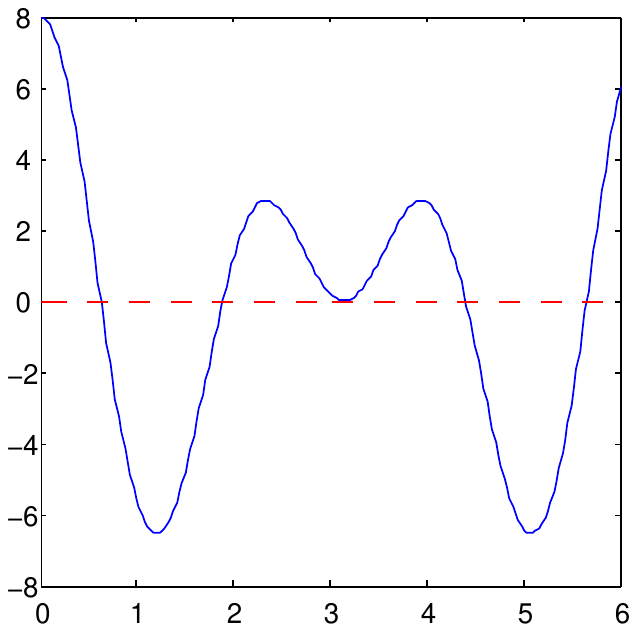} \quad\includegraphics[
natwidth=991,
natheight=750,
width=3.3in,
height=2.5in]{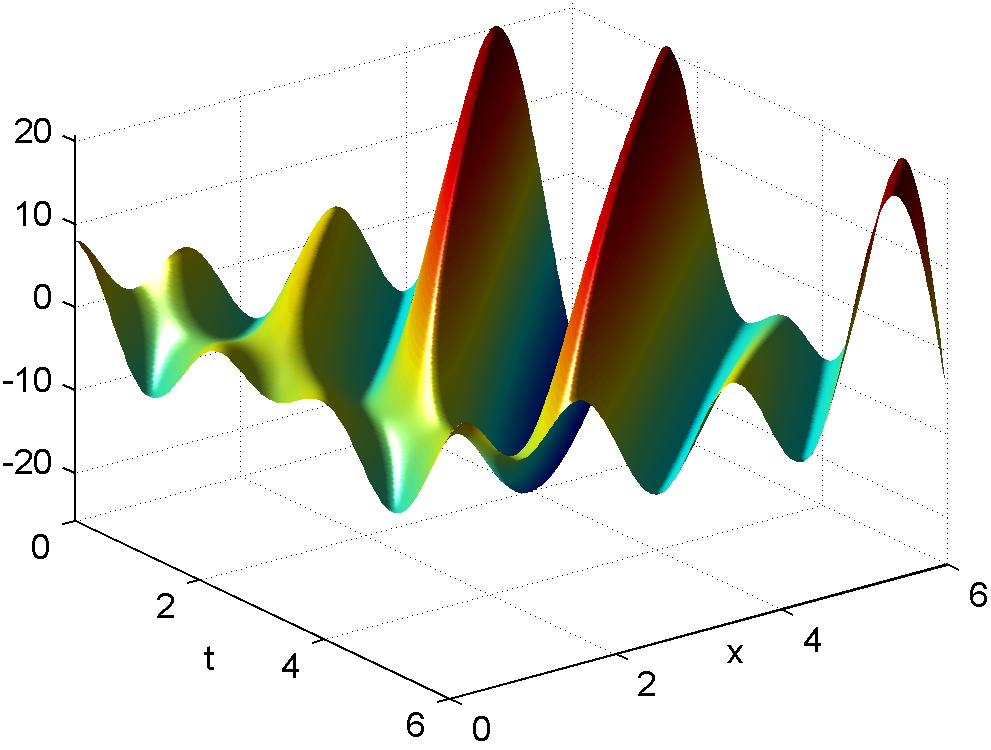}\caption{Graphs of the initial data $W^{\pm}_0(t)$ (on the left, real part in solid blue line, imaginary part in dashed red line) and of
$\mathcal{E}(x,t)$ (on the right) from
Example \ref{Ex2}.}%
\label{fig:Ex2}%
\end{figure}
\end{example}

\begin{example}\label{Ex3}
For this example let us consider the system \eqref{Maxwell1+1} with the permittivity of the form $\varepsilon(x)=(\alpha x+\beta)^{-2+2\ell}$ where $\ell\ne 0$ and $\alpha$, $\beta\in\mathbb{R}$ are such that $\alpha x+\beta>0$ on the interval of interest. Then
$\xi(x)=\sqrt\mu\int_0^x \sqrt{\varepsilon(s)}\,ds = \frac{\sqrt\mu}{\alpha\ell}\big((\alpha x+\beta)^{\ell}-\beta^{\ell}\big)$. Hence
\[
x=\frac{\beta}{\alpha}\left(\left(1+\frac{\alpha \ell \xi}{\sqrt\mu \beta^{\ell}}\right)^{1/\ell}-1\right),\quad \widetilde c(\xi)=\frac 1{\sqrt\mu} \left(\frac{\alpha\ell \xi}{\sqrt\mu } +\beta^{\ell}\right)^{(1-\ell)/\ell}\quad \text{and}\quad \widetilde f(\xi)=\frac 1{\bigl(1+\frac{\alpha \ell }{\sqrt\mu \beta^{\ell}}\xi\bigr)^{\frac{1-\ell}{2\ell}}}.
\]

In \cite{KrT2012} we show how one can construct the transmutation operators $T_f$ and $T_{1/f}$ when $f=(1+\xi)^n$, $n\in\mathbb{Z}$. The procedure can be easily generalized for functions of the form $f(\xi)=(1+c\xi)^n$, $n\in\mathbb{Z}$. Hence for values of $\ell$ of the form $\ell = \frac 1{1-2n}$, $n\in\mathbb{Z}$ one can explicitly construct the pair of transmutation operators $T_f$ and $T_{1/f}$ and obtain the solution of \eqref{Solution initial problem}.

For the numerical experiment we took $\ell = 1/5$, $\alpha = 5$, $\beta = 1$ and $\mu = 1$. For such parameters the function $\widetilde f$ is equal to $\frac 1{(1+\xi)^2}$ and the integral kernels of the transmutation operators $T_f$ and $T_{1/f}$ are given by (see \cite{KrT2012})
\begin{equation*}
    \mathbf{K}_{\widetilde f}(\xi, t) = \frac{(3t-1)(\xi+1)^2 - 3(t-1)^2 (t+1)}{4(\xi+1)^2}\quad \text{and}\quad \mathbf{K}_{1/\widetilde f}(\xi, t) = \frac{3 \xi^2 +6\xi + 4 -3t^2 + 2t}{4(\xi+1)}.
\end{equation*}
We considered the interval $[0,2]$ for $x$ and the interval $[-2,2]$ for $t$ and took the Gaussian pulse $W_0^+(t)=0$, $W_0^-(t)=e^{-4t^2}$ as the initial condition. For such initial condition the expression \eqref{Solution initial problem} can be evaluated in the terms of $\texttt{erf}$ function. The approximate solution was computed using \eqref{TfWfin} and \eqref{TfWmfin}. We used 2001 points to represent the permittivity $\varepsilon(x)$. The formal powers were computed as it was explained in Example \ref{Ex1} while to evaluate the indefinite integrals $\int z^k W_0^{\pm}(z)\,dz$, $k=0,\ldots,N$ we approximated the integrands as splines and used the function $\texttt{fnint}$ from Matlab. The main reason for such choice is that despite the uniform mesh was taken for both $x$ and $t$, the resulting mesh for $\xi$ may not be uniform leading to rather large set of values which can be taken by $t-\xi$ and $t+\xi$. All computations were performed in machine precision in Matlab 2012. On Figure \ref{fig:Ex3} we show the obtained graphs of $\mathcal{E}(x,t)$ and $\mathcal{H}(x,t)$. The absolute errors of the computed solutions were less than $1.1\cdot 10^{-14}$ and $2.6\cdot 10^{-14}$, respectively.
\begin{figure}[tbh]
\centering
\includegraphics[
width=3in,
height=2.5in,
natwidth=900,
natheight=750,
]{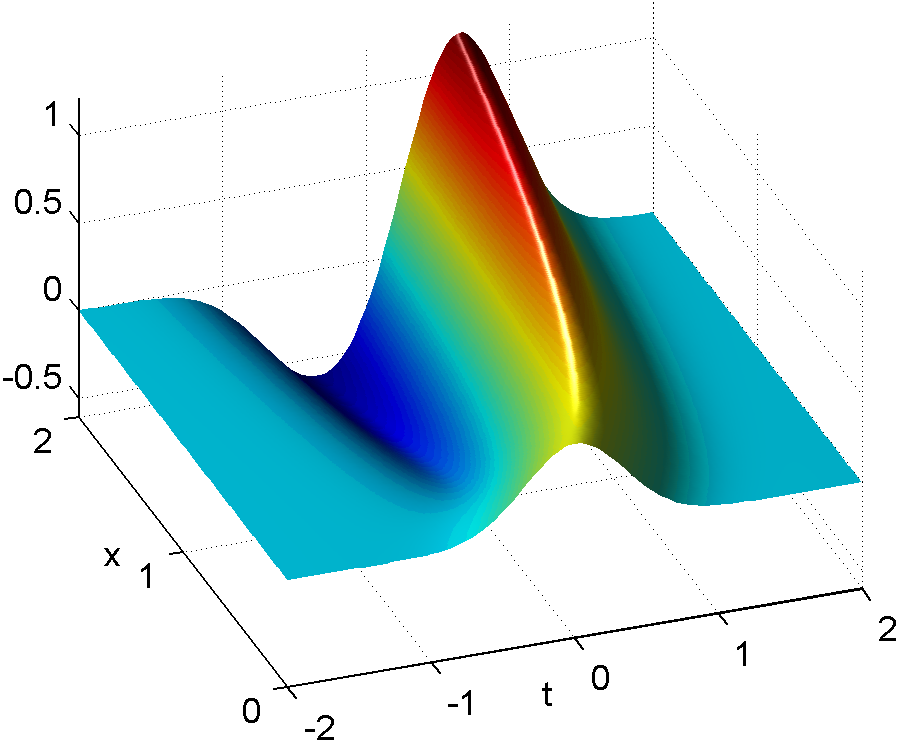} \quad\includegraphics[
natwidth=900,
natheight=750,
width=3in,
height=2.5in]{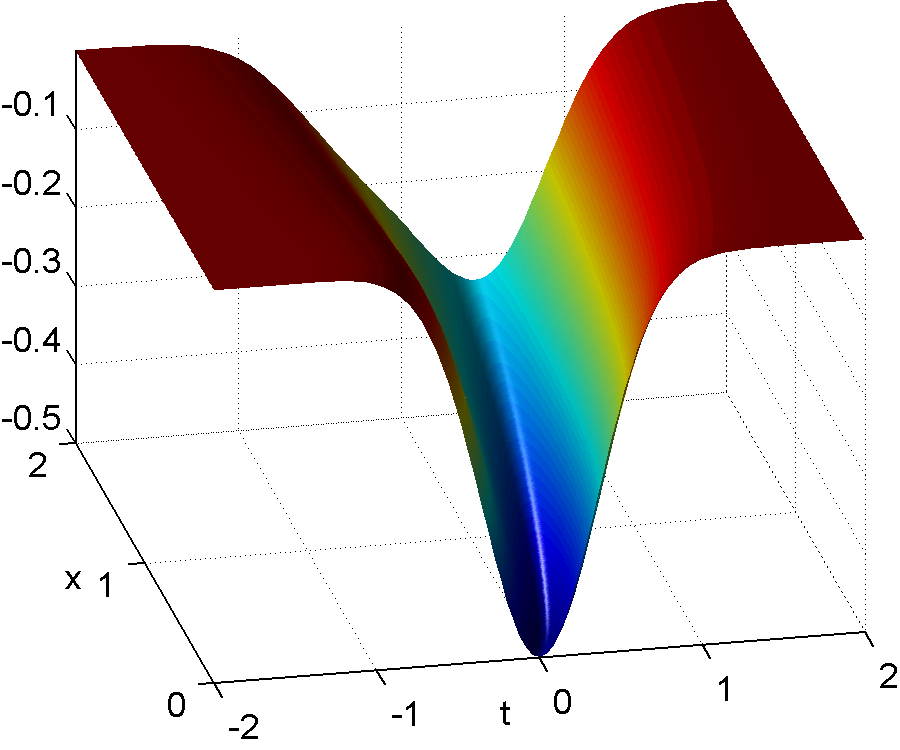}\caption{Graphs of $\mathcal{E}(x,t)$ (on the left)
and of $\mathcal{H}(x,t)$ (on the right) from
Example \ref{Ex3}.}%
\label{fig:Ex3}%
\end{figure}
\end{example}

\section{Conclusions}\label{Sect7}

A method for solving the problem of electromagnetic wave propagation through
an inhomogeneous medium is developed. It is based on a simple transformation
of the Maxwell system into a hyperbolic Vekua equation and on the solution of
this equation by means of approximate transmutation operators. In spite of
elaborate mathematical results which are behind of the proposed method, the
final representations for approximate solutions of the electromagnetic problem
have a sufficiently simple form, their numerical implementation is
straightforward and can use standard routines of such packages as Matlab.

\end{document}